# Some remarks about q-Narayana polynomials for q=−1

Johann Cigler

**Abstract**

We obtain some properties of the $q-$Narayana polynomials for $q=-1$ and compare them to corresponding properties for $q=1$.

## 1. Introduction

Let

(1) $$C_n(t;q)=\sum_{k=0}^{n-1} q^{k^2+k}\begin{bmatrix} n \\ k \end{bmatrix}\begin{bmatrix} n-1 \\ k \end{bmatrix}\frac{t^k}{[k+1]}=\sum_{k=0}^{n-1} N_{n,k}(q)t^k$$

for $n>0$ and $C_0(t;q)=1$ be $q-$Narayana polynomials. It is well known (cf. e.g. [4]) that $N_{n,k}(q)$ can be interpreted as the major index generating function of the set of Dyck paths of semi-length $n$ with $k$ valleys and that

(2) $$C_n(1;q)=\sum_{k=0}^{n-1} q^{k^2+k}\begin{bmatrix} n \\ k \end{bmatrix}\begin{bmatrix} n-1 \\ k \end{bmatrix}\frac{1}{[k+1]}=\frac{1}{[n+1]}\begin{bmatrix} 2n \\ n \end{bmatrix}$$

gives a $q-$analogue of the Catalan numbers $C_n=\dfrac{1}{n+1}\dbinom{2n}{n}$.

Here $\begin{bmatrix} n \\ k \end{bmatrix}=\begin{bmatrix} n \\ k \end{bmatrix}_q$ denotes a $q-$binomial coefficient.

For $q=-1$ identity (1) reduces to

(3) $$c_n(t)=C_n(t;-1)=\sum_{k=0}^{n} v(n,k)t^k$$

with

(4) $$v(n,k)=N_{n,k}(-1)=\binom{\left\lfloor \frac{n-1}{2}\right\rfloor}{\left\lfloor \frac{k}{2}\right\rfloor}\binom{\left\lfloor \frac{n}{2}\right\rfloor}{\left\lfloor \frac{k+1}{2}\right\rfloor}$$

and identity (2) becomes

(5) $$c_n(1)=C_n(1;-1)=\sum_{k=0}^{n} v(n,k)=\binom{n}{\left\lfloor \frac{n}{2}\right\rfloor}.$$

In [1] it has been shown that $v(n,k)$ can be interpreted as the number of symmetric Dyck paths of semi-length $n$ with $k$ valleys. In the present note we obtain some algebraic properties of the polynomials $c_n(t)$ and compare them with corresponding properties of the well-known Narayana polynomials

(6) $$C_n(t)=C(t;1)=\sum_{k=0}^{n-1}\binom{n}{k}\binom{n-1}{k}\frac{t^k}{k+1}=\sum_{k=0}^{n-1}N_{n,k}t^k.$$

The first terms of the sequence $c_n(t)$ are

$$\left(c_n(t)\right)_{n\geq 0}=\left(1,1,1+t,1+t+t^2,1+2t+2t^2+t^3,1+2t+4t^2+2t^3+t^4,\cdots\right).$$

The first Narayana polynomials are

$$\left(C_n(t)\right)_{n\geq 0}=\left(1,1,1+t,1+3t+t^2,1+6t+6t^2+t^3,1+10t+20t^2+10t^3+t^4,\cdots\right).$$

It is well known that the Hankel determinants of the Narayana polynomials satisfy

(7) $$\det\left(C_{i+j}(t)\right)_{i,j=0}^{n-1}=\det\left(C_{i+j+1}(t)\right)_{i,j=0}^{n-1}=t^{\binom{n}{2}}.$$

For the sequence $\left(c_n(t)\right)_{n\geq 0}$ we show that

(8) $$\det\left(c_{i+j}(t)\right)_{i,j=0}^{n-1}=t^{\binom{n}{2}},\quad \det\left(c_{i+j+1}(t)\right)_{i,j=0}^{n-1}=(-t)^{\binom{n}{2}}.$$

This uniquely determines the sequence $c_n(t)$.

## 2. Elementary properties of the polynomials $c_n(t)$

The polynomials $c_n(t)$ satisfy the following recursion for $n\geq 1$.

(9) $$\begin{aligned}c_{2n}(t)&=(1+t)c_{2n-1}(t),\\ c_{2n+1}(t)&=(1+t)c_{2n}(t)-tC_n\left(t^2\right).\end{aligned}$$

These identities are easily verified by comparing coefficients.

For the second identity note that $v(2n+1,2k)=v(2n,2k)+v(2n,2k-1)$ and $v(2n+1,2k+1)=v(2n,2k+1)+v(2n,2k)-N_{n,k}$.

For $t=-1$ we get

(10) $$c_{2n}(-1)=0,\quad c_{2n+1}(-1)=C_n.$$

Since $v(2n+1,2k)=\binom{n}{k}^2$ and $v(2n+1,2k+1)=\binom{n}{k}\binom{n}{k+1}$ we get

(11) $$c_{2n+1}(t)=\sum_{k=0}^{n}\binom{n}{k}^2 t^{2k}+\sum_{k=0}^{n}\binom{n}{k}\binom{n}{k+1}t^{2k+1}.$$

The polynomials

(12) $$W_n(t)=\sum_{k=0}^{n}\binom{n}{k}^2 t^k$$

are called Narayana polynomials of type B. Since

$\sum_{k=0}^{n}\binom{n}{k}\binom{n}{k+1}t^k=n\sum_{k=0}^{n}\binom{n}{k}\binom{n-1}{k}\frac{t^k}{k+1}=nC_n(t)$ we get

**Theorem 1**

*The polynomials* $c_{2n+1}(t)$ *are explicitly given by*

(13) $$c_{2n+1}(t)=W_n\left(t^2\right)+ntC_n\left(t^2\right).$$

## 3. Generating functions

The generating function

(14) $$C(t,z)=\sum_{n\geq 0}C_n(t)z^n$$

of the Narayana polynomials satisfies (cf. [4])

(15) $$C(t,z)=1-z(t-1)C(t,z)+tzC(t,z)^2$$

or equivalently

(16) $$\frac{1}{C(t,z)}=1+(t-1)z-tzC(t,z).$$

This gives the well-known formula

(17) $$C(t,z)=\sum_{n\geq 0}C_n(t)z^n=\frac{1+z(t-1)-\sqrt{\left(1+z(t-1)\right)^2-4tz}}{2tz}.$$

For the generating function of the shifted sequence

(18) $$G(t,z)=\sum_{n\geq 0}C_{n+1}(t)z^n=\frac{C(t,z)-1}{z}$$

(15) implies

(19) $$G(t,z)=\frac{1}{1-(1+t)z-tz^2G(t,z)}$$

or equivalently

(20) $$G(t,z)=1+(1+t)zG(t,z)+tz^2G(t,z)^2$$

which gives

(21) $$G(t,z)=\frac{1-(1+t)z-\sqrt{\left(1-(1+t)z\right)^2-4tz^2}}{2tz^2}.$$

Let

(22) $$c(t,z)=\sum_{n\geq 0}c_n(t)z^n$$

be the generating function of the polynomials $c_n(t)$.

The recursions (9) imply

(23) $$\left(1-(1+t)z\right)c(t,z)=1-ztC\left(t^2,z^2\right).$$

On the other hand, computations suggest that

$$\frac{1}{c(t,z)}=1-(1+t)z+tzc(-t,-z)$$

which gives

(24) $$\left(1-(1+t)z\right)c(t,z)=1-tzc\left(-t,-z\right)c(t,z).$$

Comparing (23) and (24) led to

**Theorem 2**

(25) $$c(t,z)c(-t,-z)=\sum_{n\geq 0}C_n\left(t^2\right)z^{2n}=C\left(t^2,z^2\right).$$

**Proof**

By (23) identity (25) can be written as $\frac{1-tzC\left(t^2,z^2\right)}{1-(1+t)z}\frac{1-tzC\left(t^2,z^2\right)}{1+(1-t)z}=C\left(t^2,z^2\right)$

or equivalently $\left(1-tzC\left(t^2,z^2\right)\right)^2=C\left(t^2,z^2\right)\left((1-tz)^2-z^2\right).$

Expanding this formula gives $C\left(t^2,z^2\right)=1+\left(1-t^2\right)z^2C\left(t^2,z^2\right)+z^2t^2C\left(t^2,z^2\right)^2$

which is true because of (15).

In the same way as above we get for the generating function of the shifted polynomials

(26) $$g(t,z)=\sum_{n\geq 0}c_{n+1}(t)z^n$$

the identity

(27) $$g(t,z)=1+(1+t)zg(t,z)-tz^2G\left(t^2,z^2\right).$$

**Theorem 3**

(28) $$g(t,z)g(t,-z)=G\left(t^2,z^2\right).$$

**Proof**

Identity (27) can be written as $\left(1-(1+t)z\right)g(t,z)=1-tz^2G\left(t^2,z^2\right)$. Therefore, identity (28) can be written as $\left(1-tz^2G\left(t^2,z^2\right)\right)^2=\left(1-(1+t)^2z^2\right)G\left(t^2,z^2\right)$ which can be expanded to

$G\left(t^2,z^2\right)=1+\left(1+t^2\right)z^2G\left(t^2,z^2\right)+t^2z^4G,\left(t^2,z^2\right)$ which is true because of (20).

For example, for $t=\pm1$ we get $g(1,z)=\sum_{n\geq 0}\binom{n+1}{\left\lfloor \frac{n+1}{2}\right\rfloor}z^n$ and $g(-1,z)=C(1,z^2)$.

## 4. Hankel determinants

Let $\left(a(n)\right)_{n\geq 0}$ be a sequence of real numbers with $a(0)=1$ and define a linear functional $L$ on $\mathbb{R}[x]$ by $L\left(x^n\right)=a(n)$. Then the following result is well known (cf. e.g. [2],[3])

If $\sum_{n\geq 0}a(n)z^n$ has the continued fraction expansion

(29) $$\sum_{n\geq 0}a\left(n\right)z^n=\cfrac{1}{1-s_0z-\cfrac{t_0z^2}{1-s_1z-\cfrac{t_1z^2}{1-\ddots}}}$$

where all $t_n\neq 0$ then the Hankel determinants $d_n=\det\left(a(i+j)\right)_{i,j=0}^{n-1}$ satisfy

(30) $$d_n=\prod_{j=1}^{n-1}\prod_{k=0}^{j-1}t_k.$$

By (27) and (28) $g(t,z)=\sum_{n\geq 0}c_{n+1}(t)z^n$ satisfies

(31) $$g(t,z)=\frac{1}{1-(1+t)z+tz^2g(t,-z)}=\cfrac{1}{1-(1+t)z+\cfrac{tz^2}{1+(1+t)z+\cfrac{tz^2}{\ddots}}}.$$

Here we have $s_n=(-1)^n(1+t)$ and $t_n=-t$.

This shows that their Hankel determinants are

(32) $$\det\left(c_{i+j+1}(t)\right)_{i,j=0}^{n-1}=(-t)^{\binom{n}{2}}.$$

Let us finally consider the generating function $c(t,z)$. By (24) we get

$$c(t,z)=\frac{1}{1-(1+t)z+tzc(-t,-z)}=\frac{1}{1-(1+t)z+tz\left(1-zg(-t,-z)\right)}$$
$$=\frac{1}{1-z-tz^2g(-t,-z)}.$$

Comparing with (31) we get $s_n=(-1)^n\left(1-t\right)$, with $s_0=1$ and $t_n=t$.

This implies that

$$\det\left(c_{i+j}(t)\right)_{i,j=0}^{n-1}=t^{\binom{n}{2}}. \tag{33}$$